\documentclass[12pt,reqno]{amsart}

\usepackage{amsmath, amssymb, amsfonts, amsbsy, amsthm, latexsym, epsfig, bbm, amscd, graphicx, enumerate}
\usepackage{comment}
\usepackage{xcolor}
\usepackage{dsfont}
\usepackage[mathscr]{euscript}

\allowdisplaybreaks[4]
\theoremstyle{definition}
\newtheorem{theorem}{Theorem} [section]

\newtheorem{example}[theorem]{Example}
\newtheorem{corollary}[theorem]{Corollary}
\newtheorem{lemma}[theorem]{Lemma}
\newtheorem{proposition}[theorem]{Proposition}
\newtheorem{definition}[theorem]{Definition}

\newtheorem{remark}[theorem]{Remark}

\numberwithin{equation}{section}

\newcommand{\la}{\langle}
\newcommand{\ra}{\rangle}
\newcommand{\ol}{\overline}

\newcommand{\Ker}{\mathop{\rm Ker}\nolimits}
\newcommand{\re}{\mathop{\rm Re}\nolimits}
\newcommand{\Hol}{\mathop{\rm Hol}\nolimits}
\renewcommand{\Im}{\mathop{\rm Im}\nolimits}

\def\CC{{\mathbb C}}
\def\DD{{\mathbb D}}

\def\LL{\mathcal L}

\def\RR{\mathbb R}

\def\TT{\mathbb T}
\def\beq{\begin{equation}}
\def\eeq{\end{equation}}

\def\spam{\mathop{\rm span}\nolimits}

\def\beginpf{\begin{proof}}
\def\endpf{\end{proof}}

\definecolor{vbcolor}{rgb}{0.35, 0.15, 0.6}
\newcommand{\vb}[1]{\textcolor{vbcolor}{#1}}

\newcommand{\rev}{}
\newenvironment{revised}{}{}

\newcommand\new{}
\newenvironment{newtext}{}{}

\begin{document}

\title{Frames Associated with $C_0$-Semigroups}

\author{Victor Bailey, Eva A. Gallardo-Gutiérrez, Jonathan R. Partington}

\address{Victor Bailey, Department of Mathematics,
Morehouse College,
Atlanta, GA 30314 USA}
\email{victor.bailey@morehouse.edu}

\address{Eva A. Gallardo-Guti\'errez,
Departamento de An\'alisis Matem\'atico y Matem\'atica Aplicada,
Facultad de Matem\'aticas,
Universidad Complutense de
Madrid,
Plaza de Ciencias 3, 28040 Madrid,  Spain\\
and Instituto de Ciencias Matem\'aticas ICMAT,
Madrid,  Spain }
\email{eva.gallardo@mat.ucm.es}

\address{Jonathan R. Partington, School of Mathematics, University of Leeds, Leeds LS2 9JT, UK}
\email{j.r.partington@leeds.ac.uk}

\thanks{The second and third author are partially supported by Plan Nacional  I+D grant no. PID2025-169474NB-C21, Spain. The scond author is also supported by
 the Spanish Ministry of Science and Innovation, through the ``Severo Ochoa Programme for Centres of Excellence in R\&D'' CEX2023-001347.}


\date{November 2026}


\begin{abstract}
A central question in dynamical frame theory asks when the orbit $\{\Lambda^n f\}_{n\in\mathbb{Z}_+}$ of a vector $f$ under a bounded operator $\Lambda$ on a separable infinite-dimensional Hilbert space $H$ is a frame. Christensen, Hasannasab and Philipp answered it, showing in particular that such an operator must be similar to the compression of the unilateral shift on $H^2(\mathbb{T})$ to a coinvariant subspace. We prove a continuous analogue: if a $C_0$-semigroup $\{T(t)\}_{t\ge0}\subset B(H)$ generates a continuous frame $\{T(t)\varphi\}_{t\ge0}$, then the operators $T(t)$ are simultaneously similar to the compressions of the shift semigroup to a fixed coinvariant subspace. Moreover, $\{T(t)\varphi\}_{t\ge0}$ is a continuous frame if and only if $\{C^n(I-A)^{-1}\varphi\}_{n\in\mathbb{Z}_+}$ is a frame, where $A$ is the infinitesimal generator and $C=(I+A)(I-A)^{-1}$ the cogenerator; the frame bounds of the two systems differ by a factor $2$. We also derive necessary conditions for the existence of frame vectors and discuss several examples. \new{Finally, we show that the compressed shift on a model space $K_\theta$ is the cogenerator of a $C_0$-semigroup admitting frame vectors if and only if the inner function $\theta$ has a finite angular derivative at $1$; in that case we describe all frame vectors, and the boundary reproducing kernel at $1$ generates a tight continuous frame.}
\end{abstract}

\maketitle

\section{Introduction}
Hilbert space frames were first introduced by Duffin and Schaeffer \cite{SD52} in their study of nonharmonic Fourier series, and they became a central tool in harmonic analysis after the work of Daubechies, Grossmann and Meyer \cite{DGM86}. A countable family $\{f_k\}_{k\in I}$ in a separable Hilbert space $H$ is a \emph{frame} if there are constants $0<C_1\le C_2$ such that
\[
C_1\|f\|^2\le\sum_{k\in I}|\langle f,f_k\rangle|^2\le C_2\|f\|^2\qquad(f\in H).
\]
Every frame is a complete Bessel sequence, and every $f\in H$ can be written as $f=\sum_{k\in I}c_kf_k$ with $\{c_k\}\in\ell^2(I)$. Unlike the coefficients in a basis expansion, these coefficients need not be unique. This redundancy is what makes frames flexible and robust, and it explains their wide use in signal processing and sampling theory; see \cite{Chr16,Hei11}.

Given an operator $\Lambda \in B(H)$, one often seeks conditions for $\Lambda$ to have a cyclic vector, that is, a vector $\varphi \in H$ such that $\{\Lambda^n \varphi\}_{n \in \mathbb{Z}_+}$ is complete in $H$. If, in addition, $\{\Lambda^n \varphi\}_{n \in \mathbb{Z}_+}$ is a frame (we then call $\varphi$ a frame vector for $\Lambda$), one obtains not only completeness but also basis-like expansions of all elements of $H$. Frames obtained by unilateral iterations of a bounded operator on a single vector are a type of so-called dynamical frames. A central problem in dynamical frame theory \cite{ACMT17} asks for conditions under which systems of the form $\{\Lambda^n \varphi\}_{n \in \mathbb{Z}_+}$ are frames for a separable infinite-dimensional Hilbert space. A necessary and sufficient condition was given in \cite{CMP20}; it exhibits the connection between frames and operator theory, since the bounded operators $\Lambda$ that generate a frame for $H$ must be similar to the compression of the unilateral shift on $H^2(\mathbb{T})$ to a coinvariant subspace.

More generally, a dynamical frame can be viewed as a frame representation $\{\pi(s)\varphi\}_{s \in S}$ of a countable, unital, left-cancellative semigroup $S$, where $\pi: S \to B(H)$ is a semigroup representation; for the cyclic semigroup $\mathbb{Z}_+$ these are exactly the systems $\{\Lambda^n \varphi\}_{n \in \mathbb{Z}_+}$. A characterization of such semigroup frame representations was obtained in \cite{BHKLL25}: the representation $\pi$ must be equivalent to a compression of the left regular representation on $\ell^2(S)$ to a coinvariant subspace.

The aforementioned results demonstrate the clear connection between frame theory and operator theory as the properties of certain invariant subspaces and their structure provide key tools for establishing these characterizations. When considering continuous frames generated by bounded operators, we may ask whether a similar characterization for such frames can be obtained.  In recent work  by Aldroubi et al. \cite{ALP19}, the authors considered the question of generating frames via continuous powers of a single bounded operator. Their main results in this direction show that under the assumptions that $\Lambda \in B(H)$ is invertible and self-adjoint, the system $\{\Lambda^tg\}_{g \in \mathcal{G}, \, t \in [0,1]}$ is a semi-continuous frame if and only if $\{\Lambda^tg\}_{g \in \mathcal{G}, \, t \in [0,L]}$ is a semi-continuous frame for all $L \in \mathbb R^+$ and that for $\Lambda \in B(H)$ invertible and normal, the system $\{\Lambda^tg\}_{g \in \mathcal{G}, \, t \in [0,L]}$ is a semi-continuous frame for $H$ if and only if there exists a finite set $\{t_1, t_2, \ldots t_n\}$  with $0 = t_1<t_2< \ldots t_n< L$ such that $\{\Lambda^{t_i}g\}_{g \in \mathcal{G}, \, 0 \leq i \leq n}$ is a (discrete) frame for $H$. In obtaining these results the authors utilize existing theory on the discretization of continuous frames and prior work in dynamical frame theory; however, their results do not provide a clear analogue of the characterizations of discrete dynamical frames mentioned above to the continuous frame setting.

In this paper, given a $C_0$-semigroup of operators $\{T(t)\}_{t \geq 0}\subset B(H)$ with infinitesimal generator $A$, we consider continuous frames of the form $\{T(t) \varphi\}_{t \geq 0}$, that is, frame representations of the semigroup $\mathbb{R}^+$. Associated with the semigroup is its cogenerator $C = (I+A)(I-A)^{-1}$, which is a bounded operator whenever $1$ is not in the spectrum of $A$, namely $1\notin \sigma(A)$. In this setting, we will show in Section~\ref{sec:2} (Theorem~\ref{compressed shift}) that if $\{T(t)\varphi\}_{t\ge0}$ is a continuous frame, then the operators $T(t)$ are simultaneously similar to the compressions of the shift semigroup, realized as multiplication by $e^{its}$ on $H^2(\mathbb{C}^+)$, to a fixed coinvariant subspace; this is the natural analogue of the result of Christensen, Hasannasab and Philipp. We then prove (Theorem~\ref{characterization}) that $\{T(t)\varphi\}_{t \geq 0}$ is a continuous frame for $H$ if and only if $\{C^n (I-A)^{-1}\varphi \}_{n \in \mathbb{Z}_+}$ is a discrete dynamical frame for $H$. More precisely, both systems have the same analysis operator up to a unitary map and a factor $\sqrt 2$, so that frame bounds, and in particular tightness, transfer exactly from one system to the other. Thus the theory of dynamical frames generated by the cogenerator sheds light on the continuous frames given by the semigroup.

In Section~\ref{sec:3} we will characterize the frame vectors of a $C_0$-semigroup and derive necessary conditions for their existence: the semigroup must be uniformly bounded and its adjoint must be strongly stable. In Section~\ref{sec:4} we will consider the cogenerator of several concrete semigroups (hyperbolic and parabolic composition semigroups and shift semigroups) and decide in each case whether frame vectors exist.
In the examples considered there are no frame vectors, which raises the question of whether $C_0$-semigroups with frame vectors arise naturally at all. In Section~\ref{sec:5} we answer this question by taking as cogenerator the compressed shift $S_\theta$ on a model space $K_\theta=H^2(\DD)\ominus\theta H^2(\DD)$, where $\theta$ is an inner function that is not a finite Blaschke product; the generator is then the compression to $K_\theta$ of multiplication by $(z+1)/(z-1)$. Although $S_\theta$ always admits frame vectors, for instance the reproducing kernel at the origin, the corresponding semigroup admits frame vectors if and only if $\theta$ has a finite angular derivative at $1$ (Theorem~\ref{thm:model}). In that case the frame vectors are described explicitly and the boundary reproducing kernel $k_1^\theta$ generates a tight continuous frame. By Theorem~\ref{compressed shift}, every semigroup with a frame vector is similar to a compressed shift semigroup on a model space, so these examples are, up to similarity, all the semigroups that admit frame vectors; Theorem~\ref{thm:model} determines which of these compressed semigroups do so.

\subsection*{Acknowledgements}
The first named author presented part of this work at the International Workshop on Operator Theory and its Applications (IWOTA 2026), held at Universit\'e Laval, Qu\'ebec City, on 3--7 August 2026, in the session ``Dynamical Sampling''.
In the final stage of preparing the manuscript, the authors used the AI assistant Claude (Anthropic) to check computations and to draft part of the argument in Proposition~\ref{prop:bergman}.

\section{Continuous frames generated by a $C_0$-semigroup}\label{sec:2}

Continuous frames for a Hilbert space satisfy an analogous relation to that which defines discrete frames. Also, just like a discrete frame, a continuous frame provides representations for each element of $H$ in terms of the frame vectors so that the utility of a continuous frame for reconstructing $H$ is apparent.

\begin{definition} Let $H$ be a complex Hilbert space and let $\Omega$ be a measure space with positive measure $\mu$. A continuous frame for $H$ is a family of vectors \rev{$\{h_k\}_{k \in \Omega} \subset H$} such that

$(i)$ For each $ h \in H$, \rev{$k \mapsto \langle h,h_k \rangle$} is a measurable function on $\Omega$.
\medskip

$(ii)$ There exist fixed constants $ 0 < C_{1} \leq C_{2}$ such that for each $ h \in H$,  $$
    C_{1} \|h\|^2 \leq \int_{\Omega} | \langle h , h_k \rangle |^2d\mu(k) \leq C_{2} \|h\|^2 .
    $$
\rev{If only the upper inequality in (ii) holds, the family is called a (continuous) Bessel family.}
\end{definition}
Taking $\Omega$ to be a countable set with each point having unit measure, we recover discrete frames.

As in the case for discrete frames, there are important operators defined for continuous frames which are useful for studying properties of the frame.

\begin{definition} Given a continuous frame $\{h_k\}_{k \in \Omega} \subset H$, the operator
$U: L^2(\Omega, \mu) \to H$ given \rev{(in the weak sense) by
\[
Uc= \int_{\Omega}c(k)\,h_k \, d\mu(k), \qquad\text{i.e.,}\quad \la Uc, h\ra = \int_\Omega c(k)\la h_k,h\ra\,d\mu(k)\quad (h\in H),
\]}
is the synthesis operator for the frame. Its adjoint
    $U^*: H \to L^2(\Omega, \mu)$, \rev{given by $(U^*h)(k)=\langle h,h_k \rangle$, $k\in\Omega$,}
is the analysis operator for the frame.
The map
$S := UU^* \in B(H)$
is the frame operator for the frame.
\end{definition}

It is well-known that the frame operator is a bounded and invertible operator on $H$. As in the case for discrete frames, a family of vectors $\{h_k\}_{k \in \Omega} \subset H$ is a continuous frame for $H$ if and only if the synthesis operator $U$ is well-defined \rev{(and hence bounded)} and surjective if and only if the analysis operator is a well-defined injective map onto a closed subspace of $L^2(\Omega, \mu)$.
\begin{definition} Let $H$ be a Hilbert space. We say $T: \rev{[0,\infty)} \to B(H)$ is a $C_0$-semigroup if each of the following hold:

$(i)$ $T(0) = I_H$.
\medskip

$(ii)$ $T(s+t) = T(t) T(s)$ for all $s, t \geq 0$.
\medskip

$(iii)$ $\underset{t \to 0^+}{\lim} T(t)h=h$ for every $h \in H$.

\end{definition}
It is of interest to determine exactly when a system of the form $\{T(t)\varphi\}_{t \geq 0}$, where $\{T(t)\}_{t \geq 0}$ is a $C_0$-semigroup, is a continuous frame for $H$. As there is substantial work on the theory of dynamical frames, obtaining results that establish conditions for $\{T(t)\varphi\}_{t \geq 0}$ to be a continuous frame which connect to the existing theory would be advantageous for better understanding these types of operator generated frames. In particular, following the previously mentioned results in dynamical frame theory \rev{\cite{BHKLL25,CMP20}}, we seek to establish that the operators $\{T(t)\}_{t \geq 0}$ have an analogous relation to the compression of a shift to a certain subspace.

\begin{definition}\label{def:frame-rep} Let $H$ be a Hilbert space and let $\{T(t)\}_{t \geq 0}$ be a $C_0$-semigroup on $H$. We say $\{T(t)\}_{t \geq 0}$ is a $C_0$-semigroup frame representation if there exists a vector $\varphi \in H$ such that $\{T(t)\varphi\}_{t \geq 0}$ is a continuous frame for $H$, \rev{that is, for some constants $0 < C_{1} \leq C_{2}$,}
$$
    C_{1} \|h\|^2 \leq \int_0^\infty | \rev{\langle h , T(t)\varphi\rangle} |^2dt \leq C_{2} \|h\|^2 \qquad \rev{(h\in H)}.
$$
\rev{Such a vector $\varphi$ is called a frame vector for the semigroup.}
 \end{definition}

Throughout, $\mathcal{S}_\tau$ ($\tau\ge0$) denotes the right shift on $L^2(0,\infty)$,
\[
(\mathcal{S}_\tau g)(r)=g(r-\tau)\ \ (r>\tau),\qquad (\mathcal{S}_\tau g)(r)=0\ \ (0<r\le\tau).
\]
Each $\mathcal{S}_\tau$ is an isometry, and $(\mathcal{S}_\tau^*g)(r)=g(r+\tau)$.

\begin{proposition}\label{prop:intertwining}
Let $\varphi$ be a frame vector for the $C_0$-semigroup $\{T(t)\}_{t \geq 0}$, with frame bounds $C_1\le C_2$, and let $U$ be the synthesis operator of the continuous frame $\{T(t)\varphi\}_{t \geq 0}$. Then, for every $\tau \geq 0$:
\begin{enumerate}
\item[(a)] $T(\tau)U = U\mathcal{S}_\tau$;
\item[(b)] $\Ker(U) \subset L^2(0,\infty)$ is invariant under $\mathcal{S}_\tau$;
\item[(c)] $T(\tau) = U\mathcal{S}_\tau U^\dagger$, where $U^\dagger = U^*(UU^*)^{-1}$ is the pseudoinverse of $U$;
\item[(d)] $\|T(\tau)\|\le \sqrt{C_2/C_1}$; in particular, the semigroup is uniformly bounded.
\end{enumerate}
\end{proposition}

\begin{proof}
(a) For $f \in L^2(0,\infty)$ and $h\in H$ we have
\[
\begin{aligned}
\la T(\tau)Uf,h\ra&=\la Uf,T(\tau)^*h\ra = \int_0^\infty f(t)\,\la T(\tau+t)\varphi,h\ra\,dt\\
&=\int_\tau^\infty f(r-\tau)\,\la T(r)\varphi,h\ra\,dr=\la U\mathcal{S}_\tau f,h\ra .
\end{aligned}
\]
(b) follows immediately from (a). For (c), note that $UU^\dagger=UU^*(UU^*)^{-1}=I$, so that (a) gives $T(\tau)=T(\tau)UU^\dagger=U\mathcal{S}_\tau U^\dagger$.
Finally, $\|U\|^2=\|U^*\|^2\le C_2$, and for $h\in H$,
$\|U^\dagger h\|^2=\la S^{-1}h, SS^{-1}h\ra=\la S^{-1}h,h\ra\le C_1^{-1}\|h\|^2$, where $S=UU^*\ge C_1 I$ is the frame operator. Since $\|\mathcal{S}_\tau\|=1$, (d) follows from (c).
\end{proof}

\begin{definition}
The Hardy space on the upper half-plane is the function space
$$H^2(\mathbb{C}^+) = \rev{\Big\{f \in \Hol(\{z: \Im z>0\}):\ \underset{y>0}{\sup}\int_{-\infty}^{\infty}|f(x+yi)|^2dx < \infty\Big\},}$$
and \rev{for $\tau\ge0$} the operator $M_{e^{i\tau s}} \in B(H^2(\mathbb{C}^+))$ is defined by $M_{e^{i\tau s}}f(s) = e^{i\tau s}f(s)$ for all $f \in H^2(\mathbb{C}^+)$.
\end{definition}
\begin{definition}
The Fourier transform is the unitary map $\mathcal{F}: L^2 (0 , \infty) \to H^2(\mathbb{C}^+) $ \rev{(see, e.g., \cite{MR07})} where $$(\mathcal{F}f)(s) = \frac{1}{\sqrt{2\pi}} \int_0^\infty f(t)e^{its}dt.$$
\end{definition}

\begin{lemma}\label{lem:fourier}
The map $\mathcal{F}$ satisfies $\mathcal{F}\mathcal{S}_\tau = M_{e^{i\tau s}}\mathcal{F}$ \rev{for every $\tau\ge0$ and consequently maps invariant subspaces of $\mathcal{S}_\tau$ onto invariant subspaces of the multiplication operator $M_{e^{i\tau s}}$.}
\end{lemma}
\begin{proof}
     Let $\tau \geq 0$. For any $f \in L^2 (0 , \infty)$, we have
     \[
     (\mathcal{F}\mathcal{S}_\tau f)(s) = \frac{1}{\sqrt{2\pi}} \int_\tau^\infty f(t-\tau)e^{its}dt=  \frac{1}{\sqrt{2\pi}} \int_0^\infty f(r)e^{i(r+\tau)s}dr = e^{i\tau s}(\mathcal{F} f)(s),
     \]
     so that $\mathcal{F}\mathcal{S}_\tau = M_{e^{i\tau s}}\mathcal{F}$. From this it is immediate that invariant subspaces of $\mathcal{S}_\tau$ are mapped to invariant subspaces of $M_{e^{i\tau s}}$ by $\mathcal{F}$.
\end{proof}
The following version of Beurling's theorem \rev{\cite{B48,Lax59}} will allow us to obtain a result regarding the operators $\{T(t)\}_{t \geq 0}$  generating a $C_0$-semigroup frame representation similar to the one provided in \rev{\cite{CMP20}}.

\begin{theorem}[Beurling--Lax]\label{thm:beurling-lax}
    A \rev{nonzero} closed subspace of $H^2(\mathbb{C}^+)$ is invariant under the operators $M_{e^{i\tau s}}$ for all $\tau \geq 0$ if and only if it has the form $\theta H^2(\mathbb{C}^+)$ for some inner function \rev{$\theta \in H^\infty(\mathbb{C}^+)$ (i.e., $\theta$ is a bounded analytic function on $\mathbb{C}^+$ with $|\theta|=1$ almost everywhere on $\mathbb{R}$).}
\end{theorem}

\begin{lemma}\label{lem:compression}
    Let $\theta \in \rev{H^\infty(\mathbb{C}^+)}$ be inner and define $K_\theta =\rev{H^2(\mathbb{C}^+)\ominus\theta H^2(\mathbb{C}^+)}$.
    Then $P_{K_{\theta}}M_{e^{i\tau s}} = P_{K_{\theta}}M_{e^{i\tau s}}P_{K_{\theta}}$ for every $\tau \geq 0$. \rev{Consequently, $\{P_{K_\theta}M_{e^{i\tau s}}|_{K_\theta}\}_{\tau\ge0}$ is a $C_0$-semigroup on $K_\theta$.}
\end{lemma}
\begin{proof}
    Let $f \in H^2(\mathbb{C}^+)$ and let $\tau \geq 0$. Since $H^2(\mathbb{C}^+) = \theta H^2(\mathbb{C}^+) \oplus K_{\theta}$ we have $$P_{K_{\theta}}M_{e^{i\tau s}} f = P_{K_{\theta}}M_{e^{i\tau s}}P_{K_{\theta}}f +P_{K_{\theta}}M_{e^{i\tau s}}P_{\theta H^2(\mathbb{C}^+)}f  = P_{K_{\theta}}M_{e^{i\tau s}}P_{K_{\theta}}f,$$
    since $\theta H^2(\mathbb{C}^+)$ is invariant under $M_{e^{i\tau s}}$. \rev{The semigroup property of the compressions follows, since $P_{K_\theta}M_{e^{i(\tau+\sigma)s}}=P_{K_\theta}M_{e^{i\tau s}}P_{K_\theta}M_{e^{i\sigma s}}$ on $K_\theta$, and strong continuity is inherited from $\{M_{e^{i\tau s}}\}_{\tau\ge0}$.}
\end{proof}

\begin{revised}
\begin{theorem} \label{compressed shift}
Let $\varphi$ be a frame vector for the $C_0$-semigroup $\{T(t)\}_{t \geq 0}$. Then there exist an infinite-dimensional closed subspace $K\subset H^2(\mathbb{C}^+)$, invariant under $M_{e^{i\tau s}}^*$ for every $\tau\ge0$, and a bounded invertible operator $W:K\to H$ such that
\[
T(\tau) = W\,P_{K}M_{e^{i\tau s}}|_{K}\,W^{-1}\qquad\text{for every }\tau\ge0 .
\]
Moreover, $K=K_\theta$ for some inner function $\theta$, or $K=H^2(\mathbb{C}^+)$; the latter case does not occur (see Remark~\ref{rem:converse}).
\end{theorem}

\begin{proof}
Let $U$ be the synthesis operator of the continuous frame $\{T(t) \varphi\}_{t \geq 0}$ and set $V = U \mathcal{F}^*: H^2(\mathbb{C}^+)\to H$, which is bounded and surjective. Since $VV^*=UU^*$, we have
$$V^{\dagger} = V^*(VV^*)^{-1} = \mathcal{F}U^* (UU^*)^{-1} = \mathcal{F}U^{\dagger}. $$
By Proposition~\ref{prop:intertwining}(c) and Lemma~\ref{lem:fourier}, $\mathcal{S}_\tau = \mathcal{F}^*M_{e^{i\tau s}}\mathcal{F}$ and
$$T(\tau) =  U \mathcal{F}^*M_{e^{i\tau s}}\mathcal{F} U^{\dagger} = VM_{e^{i\tau s}} V^{\dagger} \qquad(\tau\ge0).$$
By Proposition~\ref{prop:intertwining}(b) and Lemma~\ref{lem:fourier}, $\Ker(V)=\mathcal{F}(\Ker U)$ is invariant under $M_{e^{i\tau s}}$ for every $\tau\ge0$. Let $K=(\Ker V)^\perp$. If $\Ker V=\{0\}$, then $K=H^2(\mathbb{C}^+)$; otherwise, by Theorem~\ref{thm:beurling-lax}, $\Ker V=\theta H^2(\mathbb{C}^+)$ for some inner function $\theta$, and $K=K_\theta$. In both cases $K$ is invariant under every $M_{e^{i\tau s}}^*$, being the orthogonal complement of an invariant subspace.

Set $W = V|_{K}$. Then $W:K\to H$ is a bounded bijection, hence invertible by the open mapping theorem, and $K$ is infinite-dimensional because $H$ is. Moreover $V=WP_K$, and since $\operatorname{Ran}V^\dagger\subset\operatorname{Ran}V^*\subset K$ and $VV^\dagger=I$, we get $V^\dagger=W^{-1}$. Therefore
$$T(\tau) = VM_{e^{i\tau s}} V^{\dagger} = WP_{K}M_{e^{i\tau s}} W^{-1}  = WP_{K}M_{e^{i\tau s}}|_{K}\, W^{-1}$$
for every $\tau \geq 0$.
\end{proof}
\end{revised}
Theorem \ref{compressed shift} above establishes that the only operators $\{T(t)\}_{t \geq 0}$ that can generate a $C_0$-semigroup frame are those \rev{which are simultaneously similar to} the compression of an operator unitarily equivalent to the shift semigroup operator $\mathcal{S}_\tau$ to some coinvariant subspace of that operator. With this result we build on the established connection between frame theory and operator theory as once again the frame generating operators for $H$ must be closely related to the compression of a shift.

\begin{revised}
\begin{remark}\label{rem:converse}
Theorem \ref{compressed shift} establishes the necessity of the relation between a $C_0$-semigroup frame representation $\{T(t)\}_{t \geq 0}$ and the compressed shifts. Conversely, assume that $T(\tau)=W P_K M_{e^{i\tau s}}|_K W^{-1}$ for all $\tau\ge0$, where $K$ is a coinvariant subspace and $W:K\to H$ is bounded and invertible. Since similarity preserves continuous frames, $\varphi\in H$ is a frame vector for $\{T(t)\}_{t\ge0}$ if and only if $W^{-1}\varphi$ is a frame vector for the compressed semigroup $\{P_K M_{e^{i\tau s}}|_K\}_{\tau\ge0}$ on $K$. Hence the converse of Theorem~\ref{compressed shift} holds precisely for those coinvariant subspaces $K$ for which the compressed semigroup admits a frame vector.

In contrast with the discrete case, where $P_K\mathbf{1}$ is always a frame vector for the compression of the unilateral shift to a coinvariant subspace $K\subset H^2(\mathbb{D})$, there is no canonical candidate here, since constant functions do not belong to $H^2(\mathbb{C}^+)$. In fact, for $K=H^2(\mathbb{C}^+)$ the compressed semigroup is $\{M_{e^{i\tau s}}\}_{\tau\ge0}$ itself, which is unitarily equivalent to the shift semigroup on $L^2(0,\infty)$ and admits no frame vectors (see Subsection~\ref{subsec:shift}). Consequently, in Theorem~\ref{compressed shift} the synthesis operator $U$ is never injective, and $K=K_\theta$ for some inner function $\theta$.
\end{remark}
\end{revised}

We now turn our attention toward applying properties of a $C_0$-semigroup in order to study discrete systems generated by operators associated with the $C_0$-semigroup. In particular, we seek to consider operators defined by the infinitesimal generator of the semigroup and determine whether it is possible to establish a clear relation between discrete dynamical frames via these operators and the continuous frames generated by the semigroup. In doing so, we obtain a characterization of $C_0$-semigroup frames in terms of dynamical frames via the cogenerator of the semigroup.

\begin{definition}
Let $\{T(t)\}_{t \geq 0}$ be a $C_0$-semigroup \rev{on} a Hilbert space $H$. The infinitesimal generator of the semigroup is the densely defined, closed operator on $H$ given by $$Ah = \underset{t \to 0^+}{\lim}\frac{1}{t}[T(t) - I]h,$$
\rev{with domain $D(A)$ consisting of those $h\in H$ for which the limit exists.}
\end{definition}

Associated with the infinitesimal generator of a semigroup is an operator known as the cogenerator of the semigroup $C = (I+A)(I-A)^{-1}$. \rev{This operator is bounded whenever $ 1\notin \sigma(A)$, since $C=-I+2(I-A)^{-1}$. (In \cite{SNF10} the cogenerator is defined as $(A+I)(A-I)^{-1}=-C$; the sign plays no role in frame questions, since $\{C^n\psi\}$ and $\{(-C)^n\psi\}$ are simultaneously frames.)}
\begin{definition}
The Hardy space on the right half-plane is the function space
$$H^2(\mathbb{C}_+) = \rev{\Big\{f \in \Hol(\{z: \re z>0\}):\ \underset{x>0}{\sup}\int_{-\infty}^{\infty}|f(x+yi)|^2dy < \infty\Big\}.}$$
\end{definition}
\begin{definition}
The Laplace transform is the \rev{bijective} map $\mathcal{L}: L^2 (0 , \infty) \to H^2(\mathbb{C}_+) $ where $$(\mathcal{L}f)(s) = \int_0^\infty f(t)e^{-st}dt.$$
\rev{By the Paley--Wiener theorem, $\|\mathcal{L}f\|_{H^2(\mathbb{C}_+)}=\sqrt{2\pi}\,\|f\|_{L^2(0,\infty)}$.}
\end{definition}

\begin{revised}
\begin{lemma} \label{bddsurjlemma}
Let $\{T(t)\}_{t \geq 0}$ be a $C_0$-semigroup and let $\varphi\in H$ be such that $\{T(t)\varphi\}_{t \geq 0}$ is a Bessel family, with synthesis operator $U$. For $F=\LL f\in H^2(\CC_+)$ write
\[
F(-A)\varphi := Uf=\int_0^\infty f(t)\,T(t)\varphi\,dt .
\]
Then $\{T(t)\varphi\}_{t \geq 0}$ is a continuous frame for $H$ if and only if the bounded map
$\Psi_\varphi:H^2(\CC_+) \to H$, $\Psi_\varphi F=F(-A)\varphi$, is surjective.
\end{lemma}

\begin{proof}
We have $\Psi_\varphi=U\LL^{-1}$, where $(2\pi)^{-1/2}\LL$ is unitary; so $\|\Psi_\varphi\|=(2\pi)^{-1/2}\|U\|$, and $\Psi_\varphi$ is surjective if and only if $U$ is. Since a Bessel family is a continuous frame if and only if its synthesis operator is surjective, the result follows.
\end{proof}
\end{revised}

The notation $F(-A)\varphi$ is justified by the Hille--Phillips functional calculus: for instance, if the semigroup is uniformly bounded and $f(t)=e^{-\lambda t}$ with $\lambda>0$, then $F(s)=(\lambda+s)^{-1}$ and $F(-A)\varphi=(\lambda-A)^{-1}\varphi$. Lemma~\ref{bddsurjlemma} thus expresses the frame property as an interpolation problem in $H^2(\CC_+)$. Transferring this problem to $H^2(\DD)$ by means of the isomorphism $F\mapsto G$,
 \beq\label{eq:cayley}
 G(z)= (1+z)^{-1} F\Big(\frac{1-z}{1+z}\Big), \qquad  F(s) = \frac{2}{1+s}\, G\Big(\frac{1-s}{1+s}\Big),
 \eeq
one is formally led to $F(-A)\varphi=2\,G(C)(I-A)^{-1}\varphi$, that is, to the discrete system $\{C^n(I-A)^{-1}\varphi\}_{n\in\mathbb{Z}_+}$. Under \eqref{eq:cayley} and $\LL^{-1}$, the monomials $z^n$ correspond to the Laguerre functions $e^{-t}L_n(2t)$, where
\[
L_n(x)=\sum_{k=0}^n\binom nk\frac{(-x)^k}{k!}
\]
are the Laguerre polynomials. The following lemma makes this idea precise; recall that the functions $e_n(t)=\sqrt2\,e^{-t}L_n(2t)$, $n\in\mathbb{Z}_+$, form an orthonormal basis of $L^2(0,\infty)$ (see, e.g., \cite[Sec.~5.7]{Sze75}).

\begin{revised}
\begin{lemma}\label{lem:laguerre}
Let $\{T(t)\}_{t \geq 0}$ be a uniformly bounded $C_0$-semigroup with generator $A$. Then $1\notin\sigma(A)$ and, for every $\varphi\in H$, setting $\psi=(I-A)^{-1}\varphi$ and $C=(I+A)(I-A)^{-1}$:
\begin{enumerate}
\item[(i)] $C^n\psi = (-1)^n\displaystyle\int_0^\infty e^{-t}L_n(2t)\,T(t)\varphi\,dt$ for every $n\in\mathbb{Z}_+$;
\item[(ii)] for every $h\in H$,
\beq\label{eq:laguerre-identity}
\sum_{n=0}^\infty |\la h, C^n\psi\ra|^2=\frac12\int_0^\infty |\la h,T(t)\varphi\ra|^2\,dt,
\eeq
where both sides may be infinite.
\end{enumerate}
\end{lemma}

\begin{proof}
Let $M=\sup_{t\ge0}\|T(t)\|$. By the Hille--Yosida theorem, every $\lambda>0$ belongs to the resolvent set of $A$, and $(\lambda-A)^{-(k+1)}=\int_0^\infty \frac{t^k}{k!}e^{-\lambda t}T(t)\,dt$ for $k\in\mathbb{Z}_+$, the integrals converging absolutely \cite[Chap.~II]{EN00}. Since $C=2(I-A)^{-1}-I$, the binomial formula gives
\[
C^n(I-A)^{-1}=\sum_{k=0}^n\binom nk(-1)^{n-k}2^k(I-A)^{-(k+1)}
=\int_0^\infty (-1)^n e^{-t}\sum_{k=0}^n\binom nk\frac{(-2t)^k}{k!}\,T(t)\,dt ,
\]
which is (i).

For (ii), fix $h\in H$ and let $g(t)=\la h,T(t)\varphi\ra$, a continuous function with $|g|\le M\|h\|\|\varphi\|$. Since the Laguerre functions are real-valued, (i) yields
\[
\la h, C^n\psi\ra=\frac{(-1)^n}{\sqrt2}\int_0^\infty g(t)\,e_n(t)\,dt=:\frac{(-1)^n}{\sqrt2}\,c_n .
\]
It remains to show that $\sum_n|c_n|^2=\|g\|_{L^2(0,\infty)}^2$ in $[0,\infty]$. If $g\in L^2(0,\infty)$, this is Parseval's identity. Conversely, assume $\sum_n|c_n|^2<\infty$, let $k=\sum_n c_ne_n\in L^2(0,\infty)$ and $d=g-k$. Then $\int_0^\infty d\,e_n\,dt=0$ for every $n$, and since $\spam\{e_0,\dots,e_N\}=\spam\{e^{-t}t^m: 0\le m\le N\}$, also $\int_0^\infty t^me^{-t}d(t)\,dt=0$ for every $m\in\mathbb{Z}_+$. The Laplace transform $D(s)=\int_0^\infty d(t)e^{-st}\,dt$ is analytic in $\re s>0$ (as $g$ is bounded and $k\in L^2$), and $D^{(m)}(1)=(-1)^m\int_0^\infty t^m e^{-t}d(t)\,dt=0$ for all $m$. Hence $D\equiv0$, so $d=0$ a.e. by the uniqueness theorem for the Laplace transform, and $g=k\in L^2(0,\infty)$. This proves \eqref{eq:laguerre-identity}.
\end{proof}

\begin{theorem} \label{characterization}
Let $\{T(t)\}_{t \ge 0}$ be a $C_0$-semigroup and let $A$ be its infinitesimal generator. Assume $ 1 \notin \sigma(A)$, let $C = (I+A)(I-A)^{-1}$ be its cogenerator, and for $\varphi\in H$ let $\psi = (I-A)^{-1}\varphi$. Then $\{T(t)\varphi\}_{t \ge 0}$ is a continuous frame for $H$ if and only if $\{C^n \psi \}_{n \in \mathbb{Z}_+}$ is a frame for $H$. Moreover, if one of these systems is a frame, then the semigroup is uniformly bounded and
\[
\sum_{n=0}^\infty |\la h, C^n\psi\ra|^2=\frac12\int_0^\infty |\la h,T(t)\varphi\ra|^2\,dt \qquad (h\in H).
\]
Consequently, $C_1$ and $C_2$ are frame bounds for $\{T(t)\varphi\}_{t\ge0}$ if and only if $C_1/2$ and $C_2/2$ are frame bounds for $\{C^n\psi\}_{n\in\mathbb{Z}_+}$. In particular,  one system is tight if and only if the other is.
\end{theorem}

\begin{proof}
By Lemma~\ref{lem:laguerre}, it suffices to show that the semigroup is uniformly bounded whenever one of the two systems is a frame.

If $\{T(t)\varphi\}_{t \ge 0}$ is a continuous frame, this is Proposition~\ref{prop:intertwining}(d).

Assume now that $\{C^n\psi\}_{n\in\mathbb{Z}_+}$ is a frame, and let $\Phi:\ell^2(\mathbb{Z}_+)\to H$ be its synthesis operator, which is bounded and surjective. If $R$ denotes the unilateral shift on $\ell^2(\mathbb{Z}_+)$, then $C\Phi=\Phi R$, so $\Ker\Phi$ is invariant under $R$. Arguing as in the proof of Theorem~\ref{compressed shift} (see also \cite{CMP20}), with $K=(\Ker\Phi)^\perp$ and $\Phi_0=\Phi|_K$, we get $C=\Phi_0X\Phi_0^{-1}$, where $X=P_KR|_K$ is a contraction. Now $-1$ is not an eigenvalue of $C$, because $Cx=-x$ means $2(I-A)^{-1}x=0$. Hence $Y:=-X$ is a contraction which does not have $1$ as an eigenvalue, and by \cite[Chap.~III, Sec.~8]{SNF10} the operator $B=(Y+I)(Y-I)^{-1}=(X-I)(X+I)^{-1}$ generates a $C_0$-semigroup of contractions. On the other hand, $C+I=2(I-A)^{-1}$ and $C-I=2A(I-A)^{-1}$, so $A=(C-I)(C+I)^{-1}$ on $D(A)=\operatorname{Ran}(C+I)$; consequently $\Phi_0^{-1}A\Phi_0=B$. Since a generator determines its semigroup, $T(t)=\Phi_0e^{tB}\Phi_0^{-1}$, and $\|T(t)\|\le\|\Phi_0\|\,\|\Phi_0^{-1}\|$ for all $t\ge0$.
\end{proof}
\end{revised}

\begin{revised}
The same argument applied to $\{T(t)\}_{t\ge0}$ and $\{T(-t)\}_{t\ge0}$ yields a version for groups, which is relevant for the examples in Section~\ref{sec:4}.

\begin{corollary}\label{cor:groups}
Let $\{T(t)\}_{t\in\RR}$ be a uniformly bounded $C_0$-group with generator $A$. Then $\pm1\notin\sigma(A)$, the cogenerator $C=(I+A)(I-A)^{-1}$ is invertible and, for $\varphi\in H$ and $\psi=(I-A)^{-1}\varphi$,
\[
\sum_{n\in\mathbb{Z}} |\la h, C^n\psi\ra|^2=\frac12\int_{\RR} |\la h,T(t)\varphi\ra|^2\,dt \qquad (h\in H).
\]
In particular, $\{T(t)\varphi\}_{t\in\RR}$ is a continuous frame for $H$ (with respect to Lebesgue measure on $\RR$) if and only if $\{C^n\psi\}_{n\in\mathbb{Z}}$ is a frame for $H$, and the frame bounds correspond as in Theorem~\ref{characterization}.
\end{corollary}

\begin{proof}
Both $A$ and $-A$ generate uniformly bounded semigroups, so $\sigma(A)\subset i\RR$. Lemma~\ref{lem:laguerre} applied to $\{T(t)\}_{t\ge0}$ gives the sum over $n\ge0$ and the integral over $(0,\infty)$. The semigroup $\{T(-t)\}_{t\ge0}$ has generator $-A$ and cogenerator $(I-A)(I+A)^{-1}=C^{-1}$, and $(I+A)^{-1}\varphi=C^{-1}\psi$. So Lemma~\ref{lem:laguerre} gives $\sum_{n\ge0}|\la h,C^{-n-1}\psi\ra|^2=\frac12\int_{-\infty}^0|\la h,T(t)\varphi\ra|^2dt$. Adding both identities proves the claim.
\end{proof}

Theorem~\ref{characterization} establishes a precise connection between the continuous frames generated by a $C_0$-semigroup $\{T(t)\}_{t \ge 0}$ and the discrete dynamical frames generated by its cogenerator: the analysis operator of $\{C^n\psi\}_{n\in\mathbb{Z}_+}$ is $2^{-1/2}$ times the composition of the analysis operator of $\{T(t)\varphi\}_{t\ge0}$ with the unitary map $L^2(0,\infty)\to\ell^2(\mathbb{Z}_+)$ given by Laguerre coefficients (up to signs). Hence properties such as tightness or the optimal frame bounds can be read off from either system. Note also that the frame vectors of the semigroup are exactly the vectors $(I-A)\psi$, where $\psi\in D(A)$ is a frame vector of the cogenerator; a frame vector of $C$ need not lie in $D(A)$ (see Section~\ref{sec:4}).
\end{revised}

\begin{revised}
The final result for this section concerns discrete systems obtained by iterating other operators defined in terms of the generator. For semigroups of contractions and $\lambda>1$, iterations of the resolvent $(\lambda - A)^{-1}$ on any unit vector always satisfy the upper frame bound; as shown below, however, they never yield a frame.
\begin{proposition}\label{prop:bessel}
    Let  $\{T(t)\}_{t \ge 0}$ be a $C_0$-semigroup of contractions on a Hilbert space $H$,
with infinitesimal generator $A$. Let $\lambda>1$. Then for any unit vector $u$, the sequence \[\{(\lambda - A)^{-(n+1)}u\}_{n \geq 0}\] is a Bessel sequence with Bessel bound $\frac{1}{\lambda^2 -1}$.
\end{proposition}

\begin{proof}
By the Hille--Yosida theorem, $\|(\lambda-A)^{-1}\|\le \lambda^{-1}$, so $\|(\lambda - A)^{-(n+1)}u\|\le\lambda^{-(n+1)}$. By the Cauchy--Schwarz inequality, for every $x\in H$,
\[
\sum_{n=0}^\infty |\la x,(\lambda - A)^{-(n+1)}u\ra|^2\le \|x\|^2\sum_{n=0}^\infty \lambda^{-2(n+1)}=\frac{\|x\|^2}{\lambda^2-1}. \qedhere
\]
\end{proof}

\begin{remark}
Since $\sum_n\|(\lambda - A)^{-(n+1)}u\|^2<\infty$, the synthesis operator of the sequence in Proposition~\ref{prop:bessel} is a Hilbert--Schmidt operator, hence compact. A compact operator cannot map onto an infinite-dimensional Hilbert space, so this sequence is never a frame for $H$. The same bound can also be obtained from the representation $\la (\lambda - A)^{-(n+1)}u,x\ra=\int_0^\infty \frac{t^n}{n!}e^{-\lambda t}\la T(t)u,x\ra\,dt$ and the generating function $\sum_n (ts)^n/(n!)^2=I_0(2\sqrt{ts})$ of the modified Bessel function.
\end{remark}
\end{revised}

\section{Frame vectors of a $C_0$-semigroup}\label{sec:3}

Throughout this section $\{T(t)\}_{t\ge0}$ is a $C_0$-semigroup on $H$ and $\varphi\in H$. Since $t\mapsto T(t)\varphi$ is continuous, each function $t\mapsto\la h,T(t)\varphi\ra$ is continuous; hence the measurability condition in the definition of a continuous frame is automatic. We say that $\varphi$ is a \emph{Bessel vector} if $\{T(t)\varphi\}_{t\ge0}$ is a Bessel family, that is, if
\[
\int_0^\infty|\la h,T(t)\varphi\ra|^2\,dt<\infty\qquad\text{for every }h\in H .
\]
In this case the analysis operator
\[
W_\varphi:H\longrightarrow L^2(0,\infty),\qquad (W_\varphi h)(t)=\la h,T(t)\varphi\ra ,
\]
is bounded by the closed graph theorem: if $h_n\to h$ and $W_\varphi h_n\to g$ in $L^2(0,\infty)$, then $W_\varphi h_n\to W_\varphi h$ pointwise, while a subsequence converges to $g$ almost everywhere, so $g=W_\varphi h$. In the notation of Section~\ref{sec:2}, $W_\varphi=U^*$, and the frame operator is $S_\varphi:=W_\varphi^*W_\varphi$, that is,
\[
S_\varphi h=\int_0^\infty \la h,T(t)\varphi\ra\,T(t)\varphi\,dt\qquad\text{(weakly)}.
\]
Recall that $\varphi$ is \emph{cyclic} for the semigroup if $\overline{\operatorname{span}}\{T(t)\varphi:\ t\ge0\}=H$.

The following characterization is in line with those provided for discrete frames in Theorems 8.13 and 8.29 of \cite{Hei11}.

\begin{theorem}[Characterization of frame vectors] \label{frame vectors}
Let $\varphi\in H$ be a Bessel vector for $\{T(t)\}_{t\ge0}$. The following are equivalent:
\begin{enumerate}
\item $\varphi$ is a frame vector, i.e., $\{T(t)\varphi\}_{t\ge0}$ is a continuous frame for $H$;
\item $W_\varphi$ is bounded below: there is $c>0$ with $\|W_\varphi h\|_{L^2(0,\infty)}\ge c\|h\|$ for all $h\in H$;
\item $W_\varphi$ is injective and has closed range;
\item the frame operator $S_\varphi$ is invertible;
\item $\varphi$ is cyclic and $\operatorname{Ran}(W_\varphi)$ is closed.
\end{enumerate}
In this case, the optimal frame bounds are $\|S_\varphi^{-1}\|^{-1}$ and $\|S_\varphi\|=\|W_\varphi\|^2$.
\end{theorem}

\begin{proof}
(1)$\Leftrightarrow$(2). Since $\varphi$ is a Bessel vector, the upper frame inequality holds with $C_2=\|W_\varphi\|^2$, and the lower one is precisely (2) with $C_1=c^2$.

(2)$\Leftrightarrow$(3). This is the standard fact that a bounded operator between Hilbert spaces is bounded below if and only if it is injective with closed range. Indeed, if (2) holds, $W_\varphi$ is clearly injective, and if $W_\varphi h_n\to g$, then $(h_n)$ is a Cauchy sequence, so $h_n\to h$ and $g=W_\varphi h$. Conversely, if (3) holds, the bounded inverse theorem applied to $W_\varphi:H\to\operatorname{Ran}(W_\varphi)$ gives (2).

(2)$\Leftrightarrow$(4). Since $\la S_\varphi h,h\ra=\|W_\varphi h\|^2$, condition (2) means $S_\varphi\ge c^2I$. A positive operator $S$ satisfies $S\ge\varepsilon I$ for some $\varepsilon>0$ if and only if it is invertible, and then the largest such $\varepsilon$ is $\|S^{-1}\|^{-1}$. Similarly, the smallest upper bound is $\sup_{\|h\|=1}\la S_\varphi h,h\ra=\|S_\varphi\|$. This also proves the last assertion.

(3)$\Leftrightarrow$(5). We have $\Ker(W_\varphi)=\{T(t)\varphi:\ t\ge0\}^\perp$, so $W_\varphi$ is injective if and only if $\varphi$ is cyclic.
\end{proof}

In particular, every frame vector is cyclic, and a cyclic Bessel vector is a frame vector if and only if the range of its analysis operator is closed. The next result gives necessary conditions on the semigroup itself; they are the continuous counterparts of the facts that an operator generating a dynamical frame is power bounded and that its adjoint is strongly stable (cf. \cite{ACMT17,CMP20}).

\begin{proposition}\label{prop:necessary}
Assume that $\{T(t)\}_{t\ge0}$ admits a frame vector $\varphi$ with frame bounds $C_1\le C_2$. Then:
\begin{enumerate}
\item[(a)] $\varphi$ is cyclic;
\item[(b)] the semigroup is uniformly bounded: $\|T(t)\|\le\sqrt{C_2/C_1}$ for all $t\ge0$;
\item[(c)] the adjoint semigroup is strongly stable: $T(t)^*h\to0$ as $t\to\infty$ for every $h\in H$.
\end{enumerate}
\end{proposition}

\begin{proof}
(a) follows from Theorem~\ref{frame vectors} and (b) is Proposition~\ref{prop:intertwining}(d). For (c), Proposition~\ref{prop:intertwining}(c) gives $T(t)^*=(U^\dagger)^*\mathcal{S}_t^*U^*$, and $\|\mathcal{S}_t^*g\|^2=\int_t^\infty|g(r)|^2dr\to0$ for every $g\in L^2(0,\infty)$.
\end{proof}

Proposition~\ref{prop:necessary} already rules out several natural semigroups: unbounded ones, such as the composition semigroups in Subsections~\ref{subsec:hyperbolic} and \ref{subsec:parabolic}, and those whose adjoints are isometric, such as unitary groups restricted to $t\ge0$ or the left shift semigroup on $L^2(0,\infty)$.

\section{The cogenerator of certain $C_0$-semigroups}\label{sec:4}
In light of the characterization provided in Theorem \ref{characterization}, for a given $C_0$-semigroup $\{T(t)\}_{t\ge0}$, it may be advantageous to compute the cogenerator of the semigroup in order to check whether $\{T(t) \varphi \}_{t\ge0}$ is a continuous frame for some $\varphi \in H$. Since we know that $\{T(t) \varphi \}_{t\ge0}$ is a continuous frame if and only if $\{C^n\psi\}_{n \in \mathbb{Z}_+}$ is a discrete frame (where $\psi = (I-A)^{-1}\varphi$), by computing the cogenerator and applying previous work in dynamical frame theory we can establish whether  $\{T(t) \varphi \}_{t\ge0}$ is a continuous frame for $H$. Before doing so, it is worth checking the necessary conditions of Proposition~\ref{prop:necessary}, which settle some cases at once.

\subsection{A hyperbolic composition semigroup}\label{subsec:hyperbolic}

On $H^2(\DD)$ consider the $C_0$-group of composition operators
\[
(T_t f)(z) = f(\phi_t(z)), \qquad
\phi_t(z)=\frac{z+\tanh(\lambda t)}{1+z\tanh(\lambda t)}, \qquad t\in\RR,\ \lambda>0,
\]
whose infinitesimal generator is $(Af)(z)=\lambda(1-z^2)f'(z)$.
Since $\|T_t\|=e^{\lambda|t|}$ (see, e.g., \cite{CM95}), Proposition~\ref{prop:necessary}(b) shows that the semigroup $\{T_t\}_{t\ge0}$ admits no frame vectors. Moreover, this growth estimate shows that $\sigma(A)$ lies in the strip $\{|\re\mu|\le\lambda\}$; conversely, every $\mu$ with $|\re \mu|<\lambda$ is an eigenvalue, with eigenfunction $\big((1+z)/(1-z)\big)^{\mu/(2\lambda)}\in H^2(\DD)$. Hence $\sigma(A)=\{|\re\mu|\le\lambda\}$, and the cogenerator is defined precisely when $0<\lambda<1$.

The Koenigs map $h(z)=\frac{1}{2\lambda}\log\frac{1+z}{1-z}$ maps $\DD$ conformally onto the strip
\[
\Sigma_a=\{u\in\CC:\ |\Im u|<a\},\qquad a=\frac{\pi}{4\lambda},
\]
and satisfies $h\circ\phi_t=h+t$. Writing $\tau=h^{-1}$, i.e.\ $\tau(u)=\tanh(\lambda u)$, the map $\mathcal{C}_\tau f=f\circ\tau$ satisfies
\[
\mathcal{C}_\tau T_t = S_t\, \mathcal{C}_\tau ,\qquad (S_tF)(u)=F(u+t),
\]
so the composition group is similar to the translation group $(S_t)_{t\in\RR}$ acting on the image space $\mathcal{C}_\tau H^2(\DD)=\{F:\ (\tau')^{1/2}F\in H^2(\Sigma_a)\}$. Here $H^2(\Sigma_a)$ is the space of holomorphic $F$ on $\Sigma_a$ with
$\|F\|^2_{H^2(\Sigma_a)}=\sup_{|y|<a}\int_\RR |F(x+iy)|^2\,dx<\infty$, and on $\partial\Sigma_a$ one has $|\tau'(x\pm ia)|=2\lambda/\cosh(2\lambda x)$. Note that this weight does not commute with translations, so $(T_t)$ is not equivalent to translation on the unweighted space $H^2(\Sigma_a)$ (on which translations are uniformly bounded, whereas $\|T_t\|=e^{\lambda|t|}$). The translation group on $H^2(\Sigma_a)$ is nevertheless the natural model, and it can be diagonalized explicitly, as we now explain.

By the Paley--Wiener theorem for strips \cite[Chap.~1]{monguzzi} (see also \cite{BK07,zen}), the boundary Fourier transform $F\mapsto\widehat F$ is an isomorphism of $H^2(\Sigma_a)$ onto a weighted $L^2$ space:
\[
\|F\|_{H^2(\Sigma_a)}^2 \asymp \int_\RR |\widehat F(\xi)|^2 \cosh(2a\xi)\,d\xi,
\]
which follows from $\widehat{F(\cdot+iy)}(\xi)=e^{-y\xi}\widehat F(\xi)$ for $|y|<a$. In this model $\widehat{S_tF}(\xi)=e^{it\xi}\widehat F(\xi)$, so the generator $B=\frac{d}{du}$ of $(S_t)$ becomes multiplication by $i\xi$ (in particular $\sigma(B)=i\RR$), and the operator $W=(B+I)(B-I)^{-1}$ (that is, minus the cogenerator) becomes multiplication by the unimodular function
\[
m(\xi)=\frac{1+i\xi}{-1+i\xi}.
\]
Orbit questions for $W$ thus reduce to the following result.

\begin{theorem}[Frame criterion for the orbit $\{m^n g\}_{n\in\mathbb Z}$]\label{thm:orbit-frame-m}
Fix $a>0$, set $w(\xi)=\cosh(2a\xi)$ and $\mathcal H=L^2(\RR,w(\xi)\,d\xi)$, and let $M_m$ be the (unitary) operator of multiplication by $m$ on $\mathcal H$. For $g\in\mathcal H$ put
\[
\rho_g(\xi):=(1+\xi^2)\,w(\xi)\,|g(\xi)|^2 .
\]
Then $\{M_m^n g\}_{n\in\mathbb Z}=\{m^n g\}_{n\in\mathbb Z}$ is a frame for $\mathcal H$ if and only if there exist constants $0<c_1\le c_2<\infty$ such that
\begin{equation}\label{eq:frame-criterion}
c_1\le \rho_g(\xi)\le c_2 \qquad\text{for a.e. }\xi\in\RR.
\end{equation}
In that case the optimal frame bounds are $\pi\,\operatorname*{ess\,inf}\rho_g$ and $\pi\,\operatorname*{ess\,sup}\rho_g$.
\end{theorem}

\begin{proof}
The key point is the identity
\begin{equation}\label{eq:key-identity}
\sum_{n\in\mathbb Z}|\langle f, m^n g\rangle_{\mathcal H}|^2
=\pi\int_\RR \rho_g(\xi)\,|f(\xi)|^2\,w(\xi)\,d\xi
\qquad (f\in\mathcal H),
\end{equation}
both sides being allowed to equal $+\infty$.
To prove it, substitute $\xi=\tan(\theta/2)$, $\theta\in(-\pi,\pi)$, so that $d\xi=\tfrac12\sec^2(\theta/2)\,d\theta=\tfrac12(1+\xi^2)\,d\theta$ and $m(\xi)=-e^{i\theta}$, since $(1+i\xi)/(1-i\xi)=e^{2i\arctan\xi}$. Then
\[
\langle f,m^ng\rangle_{\mathcal H}=\int_\RR f\,\ol g\,\ol m^{\,n}\,w\,d\xi
=(-1)^n\int_{-\pi}^{\pi}\Phi(\theta)\,e^{-in\theta}\,d\theta,
\]
where
\[
\Phi(\theta):=\tfrac12\sec^2(\theta/2)\,(f\ol g w)(\tan(\theta/2)).
\]
Here $\Phi\in L^1(-\pi,\pi)$, because $\int_{-\pi}^{\pi}|\Phi|\,d\theta=\int_\RR|f||g|\,w\,d\xi\le\|f\|_{\mathcal H}\|g\|_{\mathcal H}$. By Parseval's identity for integrable functions on $(-\pi,\pi)$, valid in $[0,\infty]$,
\[
\sum_{n\in\mathbb Z}|\langle f,m^ng\rangle_{\mathcal H}|^2=2\pi\int_{-\pi}^{\pi}|\Phi(\theta)|^2\,d\theta
=2\pi\int_\RR \tfrac12(1+\xi^2)\,|f|^2|g|^2w^2\,d\xi,
\]
which is \eqref{eq:key-identity}.

If \eqref{eq:frame-criterion} holds, then \eqref{eq:key-identity} gives $\pi c_1\|f\|^2_{\mathcal H}\le\sum_n|\langle f,m^ng\rangle|^2\le\pi c_2\|f\|_{\mathcal H}^2$, so $\{m^ng\}_{n\in\mathbb Z}$ is a frame. Conversely, suppose it is a frame with bounds $C_1\le C_2$. If the set $\{\rho_g>C_2/\pi\}$ had positive measure, it would contain a set $E$ of finite positive measure, and $f=\mathds 1_E$ in \eqref{eq:key-identity} would violate the upper frame bound; hence $\rho_g\le C_2/\pi$ a.e. In the same way the lower frame bound forces $\rho_g\ge C_1/\pi$ a.e. This proves \eqref{eq:frame-criterion}, and also that $\pi\operatorname*{ess\,inf}\rho_g$ and $\pi\operatorname*{ess\,sup}\rho_g$ are the optimal frame bounds.
\end{proof}

Since frames are preserved by bounded invertible maps, Theorem~\ref{thm:orbit-frame-m} transfers directly to the strip.

\begin{corollary}[Orbit frames for the cogenerator on the strip model]\label{cor:strip-cogenerator}
Let $(S_t)_{t\in\RR}$ be the translation group on $H^2(\Sigma_a)$, with generator $B=\frac{d}{du}$, and let $W=(B+I)(B-I)^{-1}$. For $F\in H^2(\Sigma_a)$, the orbit $\{W^nF\}_{n\in\mathbb Z}$ is a frame for $H^2(\Sigma_a)$ if and only if there exist $0<c_1\le c_2<\infty$ such that
\[
c_1\le (1+\xi^2)\cosh(2a\xi)\,|\widehat F(\xi)|^2\le c_2
\qquad\text{for a.e. }\xi\in\RR .
\]
\end{corollary}

\begin{remark}\label{rem:strip}
The translation group $(S_t)_{t\in\RR}$ is uniformly bounded on $H^2(\Sigma_a)$; in fact it is unitary for the equivalent norm $\|\widehat F\|_{L^2(\RR,\cosh(2a\xi)d\xi)}$. Since the cogenerator is $-W$ and $(I-B)^{-1}$ corresponds to multiplication by $(1-i\xi)^{-1}$, Corollary~\ref{cor:groups} and Theorem~\ref{thm:orbit-frame-m} show that $\{S_tF\}_{t\in\RR}$ is a continuous frame for $H^2(\Sigma_a)$ if and only if
\[
c_1\le \cosh(2a\xi)\,|\widehat F(\xi)|^2\le c_2\qquad\text{for a.e. }\xi\in\RR,
\]
which is impossible, as $\int_\RR|\widehat F|^2\cosh(2a\xi)\,d\xi<\infty$. Thus the translation group admits no continuous frame vectors (and neither does the semigroup $(S_t)_{t\ge0}$, by Proposition~\ref{prop:necessary}(c)), although the operator $W$ has many orbit frames by Corollary~\ref{cor:strip-cogenerator}. This is consistent with Theorem~\ref{characterization}: the vectors $F$ in Corollary~\ref{cor:strip-cogenerator} satisfy $(1+\xi^2)\cosh(2a\xi)|\widehat F|^2\asymp1$, so they do not belong to the domain of $B$.
\end{remark}

\subsection{A parabolic composition semigroup}\label{subsec:parabolic}

\begin{revised}
Consider the parabolic composition semigroup on $H^2(\DD)$ given by
\[
(T_tf)(z)=f(\phi_t(z)),\qquad \phi_t(z)= \frac{t+(1-t)z}{1+t-tz} , \qquad z \in \DD,\ t \ge 0.
\]
The infinitesimal generator is given by $f \mapsto Gf'$ where $G(z)=(1-z)^2$.
The corresponding Koenigs mapping $h$ is $h(z)=z/(1-z)$, mapping $\DD$ onto
$\{ z \in \CC: \re z > -1/2 \}$, and $h(\phi_t(z))=h(z)+t$ (see \cite{siskakis}, for instance).

This semigroup is not uniformly bounded. Indeed, for the reproducing kernels $K_w(z)=(1-\ol w z)^{-1}$ of $H^2(\DD)$ we have $T_t^*K_0=K_{\phi_t(0)}$, and $\phi_t(0)=t/(1+t)$, so
\[
\|T_t\|^2\ge\|K_{\phi_t(0)}\|^2=\frac{1}{1-|\phi_t(0)|^2}=\frac{(1+t)^2}{1+2t}\longrightarrow\infty .
\]
Hence, by Proposition~\ref{prop:necessary}(b), it admits no frame vectors.

As in Subsection~\ref{subsec:hyperbolic}, the map $\tilde h=h+1/2$ is a bijection from $\DD$ onto the right half-plane $\CC_+$ with $\tilde h(\phi_t(z))=\tilde h(z)+t$, so that composition with $\tilde h^{-1}$ intertwines $(T_t)$ with the translation semigroup $F\mapsto F(\cdot+t)$, but on a weighted Hardy space of $\CC_+$ rather than on $H^2(\CC_+)$ itself. The translation semigroup on $H^2(\CC_+)$ is a contraction semigroup, and it is the natural model to consider. Since
\[
\int_0^\infty f(x)e^{-xt} e^{-sx} \, dx = \LL f(s+t)\qquad(f\in L^2(0,\infty)),
\]
the Laplace transform shows that it is equivalent to the semigroup $(T_t)_{t \ge 0}$ on $L^2(0,\infty)$ consisting of multiplication by $e^{-xt}$. Its infinitesimal generator $A$ is multiplication by $-x$, and the cogenerator $C=(I+A)(I-A)^{-1}$ is multiplication on $L^2(0,\infty)$ by $c(x)=(1-x)/(1+x)$, $x \ge 0$. Using the transformation $t=(1-x)/(1+x)$, the frame property related to the cogenerator reduces to considering systems of the form $\{g(t)t^n\}_{n\in\mathbb{Z}_+}$ in a (weighted) $L^2(-1,1)$ space. The next result shows that such systems are never frames.

\begin{proposition}\label{prop:multiplication}
The cogenerator $C$ (multiplication by $c(x)=(1-x)/(1+x)$ on $L^2(0,\infty)$) admits no frame vectors. Consequently, by Theorem~\ref{characterization}, the multiplication semigroup $\{e^{-xt}\}_{t\ge0}$ on $L^2(0,\infty)$, and hence the translation semigroup on $H^2(\CC_+)$, admit no frame vectors either.
\end{proposition}

\begin{proof}
Let $g\in L^2(0,\infty)$ and $E=[1/3,3]$, so that $|c|\le1/2$ on $E$. Let $(f_m)_m$ be an orthonormal sequence in $L^2(0,\infty)$ of functions supported in $E$. For every $n$ we have $|\la f_m,c^ng\ra|\le2^{-n}\|g\|$, and $\la f_m,c^ng\ra\to0$ as $m\to\infty$ by Bessel's inequality. By dominated convergence, $\sum_n|\la f_m,c^ng\ra|^2\to0$ as $m\to\infty$, while $\|f_m\|=1$. Hence $\{c^ng\}_{n\in\mathbb{Z}_+}$ does not satisfy the lower frame inequality.
\end{proof}
\end{revised}

\subsection{Shift semigroups}\label{subsec:shift}

We know that the shift semigroup on $L^2(0,\infty)$ is unitarily equivalent
to the multiplication semigroup $T_t F(s)=e^{-st} F(s)$ on the Hardy space
$H^2(\CC_+)$ of the right half-plane.

Its generator is the operator $A$ satisfying $AF(s)=-sF(s)$ for $F$ in its domain, and the
cogenerator $C$ satisfies
\[
CF(s)=\frac{1-s}{1+s}F(s),
\]
which by the isomorphism \rev{\eqref{eq:cayley}} with $H^2(\DD)$ translates to the
unilateral shift, namely,
$SG(z)=z G(z)$.

This has many frame vectors $u$, namely, those functions such that $u$ and $u^{-1}$
belong to $H^\infty(\DD)$. Translating that back to the half-plane we find that the frame
vectors for the cogenerator are those $v \in H^2(\CC_+)$ such that the function $(1+s)v(s)$ is bounded above and below.

However, \rev{by Theorem~\ref{characterization}, since $\psi=(I-A)^{-1}\varphi$ corresponds to $\varphi(s)/(1+s)$,} the frame vectors for the semigroup in $L^2(0,\infty)$ would have to be functions
whose Laplace transforms are bounded above and below on $\CC_+$, and this
is not possible. We can see this directly working with the semigroup
$(T_t)$ on $H^2(\CC_+)$ as it would require the mapping
$F \mapsto F(s) u(s)$ to be bounded and surjective for some $u$ in $H^2(\CC_+)$.

\begin{revised}
This illustrates the last observation after Theorem~\ref{characterization}: frame vectors $\varphi$ for the semigroup give frame vectors $\psi=(I-A)^{-1}\varphi$ for the cogenerator, but a frame vector $\psi$ of the cogenerator gives a frame vector $(I-A)\psi$ of the semigroup only if $\psi\in D(A)$. Here no frame vector of $C$ lies in $D(A)$, since $(1+s)v(s)$ cannot be bounded below and belong to $H^2(\CC_+)$ at the same time.
\end{revised}

Note that the backward shift on $L^2(0,\infty)$ has cogenerator equivalent to the
backward shift on $H^2(\DD)$, which is not $C_{\cdot 0}$ and so has no frame vectors \rev{(this also follows from Proposition~\ref{prop:necessary}(c), as the adjoint of the backward shift semigroup consists of isometries).}

Now consider the Bergman shift, using the Laplace isomorphism between the Bergman
space $A^2(\CC_+)$ of the right half-plane and the space $L^2((0,\infty),dt/t)$, with norm given by
\[
\|f\|^2 = \int_0^\infty \frac{|f(t)|^2}{t} \, dt.
\]
The isomorphism  $A^2(\CC_+) \to A^2(\DD)$, $F \mapsto G$, is given by
 \beq
 G(z)= (1+z)^{-2} F((1-z)/(1+z)), \qquad (z \in \DD),
 \eeq
 or equivalently,
 \[
 F(s) = \frac{4}{(1+s)^2} G((1-s)/(1+s)), \qquad (s \in \CC_+).
 \]
The corresponding cogenerator is again just the shift $S$ on $A^2(\DD)$ (see \cite{DGM,GP}, for instance).

\begin{revised}
We note that the constant function $e_0(z)=1$ is not a frame vector for the Bergman shift. Indeed, recall that the functions $e_n(z)=z^n$ are orthogonal in $A^2(\DD)$ with $\|e_n\|^2=1/(n+1)$. For $f(z)=\sum_{n=0}^\infty f_n z^n\in A^2(\DD)$ we have $\la f,S^ne_0\ra=f_n/(n+1)$, so
\[
\sum_{n=0}^\infty |\langle f, S^n e_0 \rangle|^2 = \sum_{n=0}^\infty \frac{|f_n|^2}{(n+1)^2},
\qquad\text{whereas}\qquad \|f\|^2=\sum_{n=0}^\infty \frac{|f_n|^2}{n+1}.
\]
Taking $f=e_N$ shows that the lower frame inequality fails, since the ratio of the two quantities is $1/(N+1)$. Thus $\{S^ne_0\}_{n\in\mathbb{Z}_+}$ is a Bessel sequence, but not a frame. Indeed, we may say a bit more:
\end{revised}

\begin{revised}
\begin{proposition}\label{prop:bergman}
The Bergman shift $B$, $Bf(z)=zf(z)$ on $A^2(\DD)$, admits no frame vectors: there is no $g\in A^2(\DD)$ such that $\{B^ng\}_{n\in\mathbb{Z}_+}$ is a frame for $A^2(\DD)$.
\end{proposition}

\begin{proof}
We use normalized area measure, so that $\|z^n\|^2=1/(n+1)$; in particular $B$ is a contraction and $\|f\|_{A^2}\le\|f\|_{H^2}$ for $f\in H^2(\DD)$. Let $S$ denote the unilateral shift on $H^2(\DD)$.

Suppose that $\{B^ng\}_{n\in\mathbb{Z}_+}$ is a frame and let $V:H^2(\DD)\to A^2(\DD)$ be its synthesis operator, $V(z^n)=B^ng$, which is bounded and surjective. Since $BV=VS$, the kernel of $V$ is $S$-invariant. Let $K=(\Ker V)^\perp$ and $V_0=V|_K$, which is invertible. Arguing as in the proof of Theorem~\ref{compressed shift} (see also \cite{CMP20}), $B=V_0\,P_KS|_K\,V_0^{-1}$. We distinguish two cases.

\emph{Case 1: $\Ker V=\{0\}$.} Then $B=V_0SV_0^{-1}$ and, as $S$ is an isometry,
$$\|B^nx\|\ge (\|V_0\|\|V_0^{-1}\|)^{-1}\|x\|$$
for all $x$ and $n$. This contradicts $\|B^n1\|=(n+1)^{-1/2}\to0$.

\emph{Case 2: $\Ker V\neq\{0\}$.} By Beurling's theorem \cite{B48}, $\Ker V=\theta H^2(\DD)$ for some inner function $\theta$, so $K=K_\theta$ and $B=V_0S_\theta V_0^{-1}$ with $S_\theta=P_{K_\theta}S|_{K_\theta}$.

For $0<r<1$ put $\theta_r(z)=\theta(rz)=\sum_k a_kr^kz^k$. Since $\sum_k|a_k|r^k<\infty$ and $B$, $S_\theta$ are contractions, $\theta_r(B)$ and $\theta_r(S_\theta)$ are given by norm-convergent power series, and $\theta_r(B)=V_0\,\theta_r(S_\theta)\,V_0^{-1}$.

Since $\theta H^2(\DD)$ is $S$-invariant, $p(S_\theta)=P_{K_\theta}M_p|_{K_\theta}$ for every polynomial $p$, hence $\theta_r(S_\theta)x=P_{K_\theta}(\theta_rx)$ for $x\in K_\theta$. As $P_{K_\theta}(\theta x)=0$,
\[
\|\theta_r(S_\theta)x\|\le\|(\theta_r-\theta)x\|_{H^2}=\Big(\int_{\TT}|\theta_r-\theta|^2|x|^2\,dm\Big)^{1/2}\longrightarrow0\qquad(r\to1^-),
\]
by dominated convergence, because $\theta_r\to\theta$ a.e.\ on $\TT$ and $|\theta_r|,|\theta|\le1$. Hence $\theta_r(B)\to0$ strongly on $A^2(\DD)$. On the other hand, $\theta_r(B)1=\theta_r\to\theta$ in $H^2(\DD)$, and therefore in $A^2(\DD)$. Thus $\theta=0$, contradicting the fact that $\theta$ is inner.
\end{proof}

\begin{remark}
The proof uses only two properties of $B$: it is not similar to an isometry, and no nonzero $\theta\in H^\infty(\DD)$ annihilates it. Any contraction with these two properties admits no frame vectors.
\end{remark}
\end{revised}

\begin{newtext}
\section{Model space semigroups}\label{sec:5}

None of the semigroups considered in Section~\ref{sec:4} admits frame vectors. In this section we describe a natural class of semigroups that do, and we determine all their frame vectors. By Theorem~\ref{compressed shift}, every semigroup with a frame vector is similar to a compressed shift semigroup on a model space; so we start from the cogenerator side and take as cogenerator the compressed shift on a model space of $H^2(\DD)$.

Let $\theta$ be an inner function that is not a finite Blaschke product, so that
$$K_\theta=H^2(\DD)\ominus\theta H^2(\DD)$$
is infinite-dimensional, and let $P_\theta$ be the orthogonal projection of $H^2(\DD)$ onto $K_\theta$. Let $S_\theta=P_\theta S|_{K_\theta}$ be the compressed shift. Since $\theta H^2(\DD)$ is $S$-invariant, $u(S_\theta)=P_\theta M_u|_{K_\theta}$ for every $u\in H^\infty(\DD)$. The reproducing kernel of $K_\theta$ at $\lambda\in\DD$ is
\[
k_\lambda^\theta(z)=\frac{1-\ol{\theta(\lambda)}\theta(z)}{1-\ol\lambda z},
\]
and in particular $k_0^\theta=P_\theta 1$. Recall that $\theta$ has a \emph{finite angular derivative} (in the sense of Carath\'eodory) at $\zeta\in\TT$ if $\theta$ has a unimodular nontangential limit $\theta(\zeta)$ at $\zeta$ and
$$\frac{\theta(z)-\theta(\zeta)}{z-\zeta}$$
has a finite nontangential limit at $\zeta$. By the Julia--Carath\'eodory theorem, this holds if and only if
$$\liminf_{z\to\zeta}\frac{1-|\theta(z)|}{1-|z|}<\infty,$$
and then $|\theta(r\zeta)-\theta(\zeta)|=O(1-r)$ as $r\to1^-$ (see, e.g., \cite[Chap.~7]{GMR16}). For $\theta=e^{i\gamma}B\,S_\mu$, with $B$ the Blaschke product with zeros $(a_n)$ and $S_\mu$ the singular inner function with measure $\mu$, this is equivalent to
\beq\label{eq:AC}
\sum_n\frac{1-|a_n|^2}{|\zeta-a_n|^2}+\int_{\TT}\frac{2\,d\mu(\xi)}{|\zeta-\xi|^2}<\infty
\eeq
(see \cite{AC70}, for instance). In that case the boundary kernel at $\zeta \in \partial \mathbb{D}$
$$k_\zeta^\theta(z)=\frac{1-\ol{\theta(\zeta)}\theta(z)}{1-\ol\zeta z}$$
belongs to $K_\theta$. Condition \eqref{eq:AC} holds trivially if $\zeta$ lies outside the closure of the zeros of $\theta$ and the support of $\mu$, but it also allows zeros accumulating at $\zeta$, provided they do so tangentially enough. 

\begin{lemma}\label{lem:model-frames}
Let $\theta$ be an inner function which is not a finite Blaschke product.
\begin{enumerate}
\item[(a)] $\{S_\theta^nk_0^\theta\}_{n\in\mathbb{Z}_+}=\{P_\theta z^n\}_{n\in\mathbb{Z}_+}$ is a Parseval frame for $K_\theta$.
\item[(b)] A vector $\psi\in K_\theta$ is a frame vector for $S_\theta$ if and only if $\psi=u(S_\theta)k_0^\theta=P_\theta u$ for some $u\in H^\infty(\DD)$ such that $u(S_\theta)$ is invertible.
\end{enumerate}
\end{lemma}

\begin{proof}
(a) Since $\theta H^2(\DD)$ is $S$-invariant, $S_\theta^nP_\theta=P_\theta S^n$ on $H^2(\DD)$, so $S_\theta^nk_0^\theta=P_\theta z^n$. For $f\in K_\theta$ we get
$$\sum_n|\la f,P_\theta z^n\ra|^2=\sum_n|\la f,z^n\ra|^2=\|f\|^2.$$

\smallskip

(b) If $u(S_\theta)$ is invertible, then $\{S_\theta^nu(S_\theta)k_0^\theta\}=u(S_\theta)\{S_\theta^nk_0^\theta\}$ is the image of a Parseval frame under an invertible operator, hence a frame.

\smallskip

Conversely, let $\psi$ be a frame vector, and let $V:H^2(\DD)\to K_\theta$ be the synthesis operator, $Vz^n=S_\theta^n\psi$, which is bounded and surjective and satisfies $VS=S_\theta V$.

For $0<r<1$ let $\theta_r(z)=\theta(rz)$. Since $\theta_r(S_\theta)x=P_\theta(\theta_rx)$ and $P_\theta(\theta x)=0$, dominated convergence on $\TT$ gives $$\|\theta_r(S_\theta)x\|\le\|(\theta_r-\theta)x\|_{H^2}\to0$$ as
$r\to1^-$, for every $x\in K_\theta$.

If $\Ker V=\{0\}$, then $S=V^{-1}S_\theta V$ and $\theta_r=\theta_r(S)1=V^{-1}\theta_r(S_\theta)V1\to0$, which is impossible since $\theta_r\to\theta$ in $H^2(\DD)$. Hence, by Beurling's theorem, $\Ker V=\theta'H^2(\DD)$ for some inner function $\theta'$, and
$$W=V|_{K_{\theta'}}:K_{\theta'}\to K_\theta$$
is invertible with $WS_{\theta'}=S_\theta W$. Then $\theta_r(S_{\theta'})=W^{-1}\theta_r(S_\theta)W\to0$ strongly, while $$\theta_r(S_{\theta'})k_0^{\theta'}=P_{\theta'}(\theta_rk_0^{\theta'})\to P_{\theta'}(\theta k_0^{\theta'}).$$
Thus
$$\theta k_0^{\theta'}=\theta-\ol{\theta'(0)}\,\theta\theta'\in\theta'H^2(\DD),$$ so
$\theta\in\theta'H^2(\DD)$. \new{Exchanging the roles of $\theta$ and $\theta'$ (note that $\theta'_r(S_{\theta'})\to0$ strongly by the same argument, and hence $\theta'_r(S_\theta)=W\theta'_r(S_{\theta'})W^{-1}\to0$ strongly), we also get $\theta'\in\theta H^2(\DD)$. Writing $\theta=\theta'h$ and $\theta'=\theta g$ with $g,h\in H^2(\DD)$, we see that $g$ and $h$ are inner and $gh=1$, so $g$ and $h$ are unimodular constants. Hence $\Ker V=\theta H^2(\DD)$.}

\smallskip

Therefore $W=V|_{K_\theta}$ is an invertible operator on $K_\theta$ commuting with $S_\theta$. \new{By Sarason's theorem on the commutant of $S_\theta$ \cite{Sar67}, every bounded operator on $K_\theta$ commuting with $S_\theta$ is of the form $u(S_\theta)$ for some $u\in H^\infty(\DD)$; thus $W=u(S_\theta)$.}

\smallskip

\new{Finally, since $(I-P_\theta)1\in\theta H^2(\DD)=\Ker V$,
\[
\psi=V1=VP_\theta1=Wk_0^\theta=u(S_\theta)k_0^\theta=P_\theta u. \qedhere
\]}
\end{proof}

The next lemma states when $k_0^\theta$ is in the range of $I-S_\theta$.
It can be found by combining Theorem~1.2 in \cite{AC70} with Lemma 3.1 in \cite{clark72} and the discussion that follows it.

\begin{lemma}\label{lem:range}
$k_0^\theta\in\operatorname{Ran}(I-S_\theta)$ if and only if $\theta$ has a finite angular derivative at $1$. In that case $k_1^\theta\in K_\theta$ and $(I-S_\theta)k_1^\theta=k_0^\theta$.
\end{lemma}

By \cite[Chap.~III, Thm.~8.1]{SNF10}, a contraction is the cogenerator, in the sense of \cite{SNF10}, of a $C_0$-semigroup of contractions if and only if $1$ is not one of its eigenvalues. The eigenvalues of $S_\theta$ are the zeros of $\theta$ in $\DD$, so this applies to $S_\theta$. Explicitly, the operator
\beq\label{eq:A-model}
A=(S_\theta+I)(S_\theta-I)^{-1},\qquad D(A)=\operatorname{Ran}(S_\theta-I),
\eeq
which is the compression to $K_\theta$ of multiplication by $(z+1)/(z-1)$, generates a $C_0$-semigroup of contractions $\{T_\theta(t)\}_{t\ge0}$ on $K_\theta$. A direct computation shows that its cogenerator, in our convention, is $C=(I+A)(I-A)^{-1}=-S_\theta$, and that $(I-A)^{-1}=\frac12(I-S_\theta)$. If $1\notin\sigma(S_\theta)$, that is, if $\theta$ extends analytically across a neighbourhood of $1$, then $A$ is bounded and $T_\theta(t)=e^{tA}$.

\begin{theorem}\label{thm:model}
Let $\theta$ be an inner function that is not a finite Blaschke product, and let $\{T_\theta(t)\}_{t\ge0}$ be the $C_0$-semigroup on $K_\theta$ generated by \eqref{eq:A-model}.
\begin{enumerate}
\item[(a)] $\{T_\theta(t)\}_{t\ge0}$ admits a frame vector if and only if $\theta$ has a finite angular derivative at $1$.
\item[(b)] In that case, the frame vectors are exactly the vectors $\varphi=2u(S_\theta)k_1^\theta=2P_\theta(u\,k_1^\theta)$ with $u\in H^\infty(\DD)$ and $u(S_\theta)$ invertible. Moreover, $\{T_\theta(t)k_1^\theta\}_{t\ge0}$ is a tight continuous frame:
\[
\int_0^\infty|\la h,T_\theta(t)k_1^\theta\ra|^2\,dt=\frac12\|h\|^2\qquad(h\in K_\theta).
\]
\end{enumerate}
\end{theorem}

\begin{proof}
Since the semigroup consists of contractions, $1\notin\sigma(A)$, and Theorem~\ref{characterization} applies with $C=-S_\theta$. Thus $\varphi$ is a frame vector for $\{T_\theta(t)\}$ if and only if
$$\psi=(I-A)^{-1}\varphi=\frac12(I-S_\theta)\varphi$$
is a frame vector for $-S_\theta$, or equivalently for $S_\theta$, since $(-S_\theta)^n\psi=(-1)^nS_\theta^n\psi$. By Lemma~\ref{lem:model-frames}, this happens if and only if $$\frac12(I-S_\theta)\varphi=u(S_\theta)k_0^\theta$$
for some $u\in H^\infty(\DD)$ with $u(S_\theta)$ invertible.

\smallskip

If such $\varphi$ and $u$ exist, then, since $u(S_\theta)^{-1}$ commutes with $S_\theta$, we get
\[
k_0^\theta=(I-S_\theta)\big[\tfrac12u(S_\theta)^{-1}\varphi\big]\in\operatorname{Ran}(I-S_\theta).
\] By Lemma~\ref{lem:range}, $\theta$ has a finite angular derivative at $1$, and $(I-S_\theta)\varphi=2u(S_\theta)(I-S_\theta)k_1^\theta=(I-S_\theta)\,2u(S_\theta)k_1^\theta$. Since $1$ is not an eigenvalue of $S_\theta$, it follows that $\varphi=2u(S_\theta)k_1^\theta$.

\smallskip

Conversely, if $\theta$ has a finite angular derivative at $1$ and $u(S_\theta)$ is invertible, Lemma~\ref{lem:range} gives $$\frac12(I-S_\theta)\,2u(S_\theta)k_1^\theta=u(S_\theta)k_0^\theta,$$
so $2u(S_\theta)k_1^\theta$ is a frame vector. Taking $u=1$, the vector $\psi=k_0^\theta$ generates a Parseval frame by Lemma~\ref{lem:model-frames}(a). By Theorem~\ref{characterization}, $\{T_\theta(t)\,2k_1^\theta\}_{t\ge0}$ is a tight frame with bound $2$, and the identity for $k_1^\theta$ follows by scaling. This concludes the proof.
\end{proof}

\begin{example}
If $\theta=\exp\big(-a\frac{1+z}{1-z}\big)$ with $a>0$, then \eqref{eq:AC} fails at $1$ and $\{T_\theta(t)\}$ has no frame vectors. In this case, $\{T_\theta(t)\}$ is unitarily equivalent to the truncated shift on $L^2(0,a)$, which is nilpotent and whose analysis operators are Hilbert--Schmidt.
\end{example}

\begin{remark}\label{rem:model}
Let $\omega(s)=(s-i)/(s+i)$, which maps $\CC^+$ onto $\DD$ with $\omega(\infty)=1$. The unitary map $H^2(\DD)\to H^2(\CC^+)$, $G\mapsto\pi^{-1/2}(s+i)^{-1}G(\omega(s))$, intertwines $S$ with $M_\omega$, which is the cogenerator, in the sense of \cite{SNF10}, of $\{M_{e^{i\tau s}}\}_{\tau\ge0}$. It maps $K_\theta$ onto $K_\Theta$ with $\Theta=\theta\circ\omega$, and it carries $\{T_\theta(t)\}$ onto the compressed semigroup $\{P_{K_\Theta}M_{e^{its}}|_{K_\Theta}\}_{t\ge0}$. Therefore, in the setting of Theorem~\ref{compressed shift} and Remark~\ref{rem:converse}, Theorem~\ref{thm:model} shows that the compressed shift semigroup on $K_\Theta\subset H^2(\CC^+)$ admits frame vectors if and only if $\Theta\circ\omega^{-1}$ has a finite angular derivative at $1$, the image of $\infty$ under $\omega$.
\end{remark}

\end{newtext}

\end{document}